\documentstyle[11pt]{article}
\def\rlx{\relax\leavevmode}
\def\inbar{\vrule height1.5ex width.4pt depth0pt}
\def\IZ{\rlx\hbox{\small \sf Z\kern-.4em Z}}
\def\IR{\rlx\hbox{\rm I\kern-.18em R}}
\def\ID{\rlx\hbox{\rm I\kern-.18em D}}
\def\IC{\rlx\hbox{\,$\inbar\kern-.3em{\rm C}$}}
\def\one{\hbox{{1}\kern-.25em\hbox{l}}}

\begin{document}
\begin{titlepage}

\hfill{UTAS-PHYS-95-54}
\vskip 1.6in
\begin{center}
{\Large {\bf Deformed boson algebras and the quantum double construction
}}
\end{center}

\normalsize
\vskip .4in

\begin{center}
D. S. McAnally$^{\dag}$  \hspace{3pt} 
and \hspace{3pt} I. Tsohantjis$^{*}$
\par \vskip .1in \noindent
${\dag}${\it Department of Mathematics, University of Queensland}\\
{\it Brisbane, Queensland, Australia 4072}\\
${*}${\it Department of Physics, University of Tasmania}\\
{\it GPO Box 252C Hobart, Australia 7001}
\end{center}
%\footnotetext{ }
\par \vskip .3in

\begin{center}
{\Large {\bf Abstract}}\\

\noindent
The quantum double construction of a $q$-deformed boson algebra possessing a 
Hopf algebra structure is carried out explicitly. The $R$-matrix thus 
obtained is compared with the existing literature.
\end{center}

\end{titlepage}

Recently there has been an increasing interest in the deformation
of Lie (super)algebras\cite{drin,j,q,kr,tol,gr1} and their quasitriangular
Hopf algebra nature\cite{sweed}, mainly because of there wide applications in 
mathematical physics. Parallel attempts to consistently $q$-deform the boson 
algebra also appeared\cite{kur,j1,mac,bid,bid2,ch,che,dav,yan2} both 
independently and in connection with quantum group realizations, 
addressing also their possible Hopf algebra nature\cite{yan,oh}. 
The main aim of this letter is, on the one hand to point out the 
ambigious validity of an $R$-matrix obtained 
from a definition of $q$-boson algebra endowed with a Hopf algebra
structure\cite{yan}, and on the other 
to demonstrate the quantum double construction\cite{drin,ros,gr} for 
this algebra which will lead to an unambiguously valid $R$-matrix. 

The $q$-deformed boson algebras, denoted here by $L$, that have been 
considered are usually 
taken to be generated by  
$a$, $a^{\dag}$
and $N$ subject to the following commutation relations:   
\begin{eqnarray}
& [N,a] = -a, \nonumber\\
& [N,a^{\dag}] = a^{\dag},
\label{1} 
\end{eqnarray}
together with {\it one} out of the following list of additional relations:
\begin{eqnarray}
& [a,a^{\dag}] = \left[ N + I\right]_{q} - 
\left[ N \right]_{q},
\label{2} 
\end{eqnarray}
\begin{eqnarray}
aa^{\dag} - q^{-1}a^{\dag}a= q^{N},
\label{3}
\end{eqnarray}
\begin{eqnarray}
aa^{\dag} - qa^{\dag}a= q^{-N}, 
\label{4}
\end{eqnarray}
\begin{eqnarray}
a^{\dag}a = [N],\quad\mbox{and} \quad aa^{\dag} = [N+I],
\label{5}
\end{eqnarray}
where $I$ is the unit of $L$ and as usual $[x]=(q^{x}-q^{-x})/(q-q^{-1})$,
and $q$ not a root of unity. When $q=1$ we obtain the well known
defining relations of the undeformed
boson algebra. It should be mentioned that the consistency of the above 
definitions is justified as they can also be obtained from $sl_{q}(2)$ by 
contraction\cite{ch,che,ng}. Generalizations of $q$--boson defining relations,
in particular that of (\ref{3}), (\ref{4}) have also been 
studied\cite{j1,b,da,jo}.
Analysis of representations of $L$ is quite rich \cite{cgp,jo}, but the
most ususally used is the $q$--Fock representation (which has been 
shown \cite{kd} to be 
isomorphic with the usual boson Fock space by expressing the 
$q$--bosons as suitable functions of the undeformed bosons) given by:
\begin{eqnarray}
|n> = ([n]!)^{1/2}(a^{\dag})^{n}|0>, \quad N|n>=n|n>, \nonumber\\
a^{\dag}|n>= [n+1]^{1/2}|n+1>, \quad a|n> = [n]^{1/2}|n-1>
\label{fock}
\end{eqnarray}
where $n=0,1,...$.
Using this representation, one can also show \cite{ng,kd,oh} the
equivalence amongst the above definitions, which {\it does not} imply, 
though, an equivalence at the abstract algebraic level (as has been 
demonstrated in \cite{oh}). 

The most important point though concerns the Hopf algebra structure
of the deformed boson algebra. Initially Hong Yan\cite{yan} showed that when 
$L$ is defined by (\ref{1}) and (\ref{2}) 
(with $N {\rightarrow} N-1/2$, see (\ref{y}) below) $L$ is a
 Hopf algebra. Later this result was generalized in \cite{oh} where (\ref{1})
was also generalized. We shall concentrate hereafter on the Hopf algebra
$L$ as defined in \cite{yan} by (\ref{1}) and a symmetrized version of
(\ref{2}), namely
\begin{eqnarray}
& [a,a^{\dag}] = \left[ N + \frac{1}{2} \right]_{q} - 
\left[ N - \frac{1}{2} \right]_{q}.
\label{y} 
\end{eqnarray}

The structure of this letter is as follows. First we give general
information on quasitriangular Hopf algebras, and also present the model of 
\cite{yan}, focussing mainly on the claimed $R$-matrix and pointing out
some inconsistencies in its properties. Then, by demonstrating the method 
of quantum double
construction, we apply it to the Hopf algebra $L$ defined by (\ref{1}) and 
(\ref{y}) to obtained a valid $R$-matrix which can be compared with that
of \cite{yan}.

Consider an algebra $A$, say over $\IC$, with multiplication 
$m:A{\otimes}A {\rightarrow} A$ 
(i.e. $m(a{\otimes}b)= ab$, ${\forall}a$ and ${\forall}b$ ${\in}A$) 
and unit $u:\IC {\rightarrow} A$ 
(i.e. $u(1)=I$, the identity on $A$) endowed with 
a Hopf algebra structure that is, having a coproduct 
$\Delta: A{\rightarrow}A{\otimes}A$, counit $\epsilon:A{\rightarrow}\IC$ 
(which is a homomorphism) and antipode $S:A{\rightarrow}A$ 
(which is an antihomomorphism i.e. $S(ab)=S(b)S(a)$) 
with the following consistency conditions:
\begin{eqnarray}
(id{\otimes}{\Delta}){\Delta}(a) = ({\Delta}{\otimes}id){\Delta}(a)
\nonumber\\
m(id{\otimes}S){\Delta}(a) = m(S{\otimes}id){\Delta}(a) =
{\epsilon}(a)1, \nonumber\\
({\epsilon}{\otimes}id){\Delta}(a) =
(id{\otimes}{\epsilon}){\Delta}(a)=a
\label{con}
\end{eqnarray}
Note also that we have $\epsilon(I)=1$, $S(I)=I$ and 
${\epsilon}(S(a))={\epsilon}(a)$, ${\forall}a{\in}A$.
Following Sweedler\cite{sweed} we write:
\begin{eqnarray}
\Delta(a) = \sum_{(a)}a^{(1)}{\otimes}a^{(2)} \quad {\forall}a{\in} A 
\quad \mbox{and generally}, \nonumber\\
\Delta_{n}(a) = ({\Delta}{\otimes}I^{\otimes(n-1)}){\Delta}_{n-1}(a)=
\sum_{(a)}a^{(1)}{\otimes}a^{(2)}\ldots{\otimes}a^{(n+1)} 
\quad \mbox{for} \quad n{\geq}2
\end{eqnarray}
Let $T$ denote the twist
map on $A{\otimes}A$, $T(a{\otimes}b)=b{\otimes}a$ and assume that $S^{-1}$, 
the inverse of the antipode, exists. Then there exists an opposite Hopf 
algebra
structure on $A$ with coproduct and antipode $T\Delta$ and $S^{-1}$ 
respectively.
According to Drinfeld\cite{drin} a Hopf algebra $A$ is called quasitriangular if there exists an
invertible element $R$ such that
\begin{eqnarray}
R = \sum_{i}a_{i}{\otimes}b_{i}\quad {\in}A{\otimes}A \nonumber\\
T\Delta(a)R=R\Delta(a) \quad {\forall}a{\in}A
\label{da}
\end{eqnarray}
Then it can be shown that $R$ satisfies the Yang--Baxter equation
\begin{eqnarray}
R_{12}R_{13}R_{23} = R_{23}R_{13}R_{12} \quad\mbox{where} \nonumber\\
(\Delta{\otimes}I)R=R_{13}R_{23}, \quad
(I{\otimes}\Delta)R=R_{13}R_{12}.
\label{dd}
\end{eqnarray}
Turning now to the algebra $L$ given by (\ref{1}), (\ref{y}),
the coproduct $\Delta$, counit $\varepsilon$ and antipode $S$ are 
 respectively as follows:
\begin{eqnarray}
\Delta(N) = N \otimes I + I \otimes N - \frac{i \alpha}{\gamma} I \otimes 
I, \nonumber\\
\Delta(a) = \left( a \otimes q^{N/2} + i q^{-N/2} \otimes a \right) 
e^{-i \alpha / 2},\nonumber\\
\Delta(a^{\dag}) = \left( a^{\dag} \otimes q^{N/2} + i q^{-N/2} \otimes 
a^{\dag} \right) e^{-i \alpha / 2}, \nonumber\\
\varepsilon(N) = \frac{i \alpha}{\gamma}, \quad \varepsilon(a) 
= \varepsilon(a^{\dag}) = 0,\quad \varepsilon(I)=1, \nonumber\\
S(N) = - N + \frac{2 i \alpha}{\gamma}I, \nonumber\\
S(a) = - q^{-1/2} a,\nonumber\\
S(a^{\dag}) = - q^{1/2} a^{\dag}. 
\label{hop}
\end{eqnarray}
where $\alpha=2{\kappa}{\pi}+{\pi}/2$,$(\kappa{\in}{\IZ})$ and $\gamma=lnq$. 
The consistency of 
these operations can be verified by direct calculation using (\ref{1}),  
(\ref{y}) and the consistency relations (\ref{con}).
Moreover the ideal $K$ of $L$ generated by 
$C= a^{\dag} a - [ N - \frac{1}{2} ]_{q}$ is not a Hopf
ideal (as $\Delta(C){\notin}K{\otimes}L+L{\otimes}K$) and 
thus the quotient algebra isomorphic to the one generated by 
(\ref{1}) and both
$a^{\dag} a = \left[ N - \frac{1}{2} \right]_q$,
$a a^{\dag} = \left[ N + \frac{1}{2} \right]_q$, is not a Hopf algebra. 

Although the Fock space (\ref{fock}), but with $N|n>=(n+1/2)|n>$, furnishes a 
representation of $L$, we can use a more general one given by:
\begin{eqnarray}
N |n> &=& (n + c) |n> \nonumber\\
a |n> &=& ([n+c-1/2]-[c-1/2])^{1/2} |n-1> \nonumber\\
a^{\dag} |n> &=& ([n+c+1/2]-[c-1/2])^{1/2} |n+1>
\label{f2}
\end{eqnarray} 
where $c$ is a non--zero complex number. If $c = \frac{1}{2}$ we obtain the 
space just mentioned above. For $|q| \neq 1$ we can choose the states 
$|n>$ to be given by 
\begin{eqnarray}
|n> = \left(\frac{i^n (q^{1/2}+q^{-1/2})^n \Gamma^{+}_{q^{1/2}}\left(2c
+\frac{i \pi}{2 \ln q}\right)}{(q^{1/2}-q^{-1/2})^n [n]_{q^{1/2}}! 
\Gamma^{+}_{q^{1/2}}\left(n+2c+\frac{i \pi}{2 \ln q}\right)}\right)^{1/2} 
{a^{\dag}}^{n}|0>, \quad \mbox{or} \nonumber\\
|n> = \left(\frac{(-i)^n (q^{1/2}+q^{-1/2})^n \Gamma^{+}_{q^{1/2}}\left(-n-2c
-\frac{i \pi}{2 \ln q}\right)}{(q^{1/2}-q^{-1/2})^n [n]_{q^{1/2}}! 
\Gamma^{+}_{q^{1/2}}\left(-2c-\frac{i \pi}{2 \ln q}\right)}\right)^{1/2} 
{a^{\dag}}^{n}|0>, 
\end{eqnarray}
whichever is well defined given a fixed value of $c$ and $\ln q$. Note that 
at least one of the above expressions is always well--defined, and if they 
are both well--defined, they are equal. 
The symmetric $q^{1/2}$--factorial $[n]_{q^{1/2}}!$ is defined
similarly to the symmetric $q$--factorial as in \cite{mac,qexp}, and 
the symmetric $q^{1/2}$--gamma function $\Gamma^{+}_{q^{1/2}}(z)$
similarly to the symmetric $q$--gamma function $\Gamma^{+}_{q}(z)$ as in
\cite{qgam} and $<0|0> = 1$.
If we now consider the quantum algebra $U_{q^{1/2}}(sl(2))$ defined by: 
\begin{eqnarray}
{[}e , f{]}=[h]_{q^{1/2}}, \quad {[}h , e{]}=2e,
\quad {[}h , f{]}=-2f,
\label{sl}
\end{eqnarray}
we can easily verify that the following expressions of $h, \ e, \ f$
in terms of the $q$--bosons do indeed satisfy (\ref{sl}):
\begin{eqnarray}
h=2 N - \frac{2 i \alpha}{\gamma}, \quad e=\lambda a^{\dag},\quad 
f=\frac{i (q^{1/2} + q^{-1/2})}{\lambda (q^{1/2} - q^{-1/2})} a, 
\end{eqnarray}
where $\lambda$ is some constant. In this realization of
$U_{q^{1/2}}(sl(2))$, the central element $C$ is the deformed 
quadratic Casimir.
Conversely, given $U_q(sl(2))$ as in (\ref{sl})(with $q^{1/2}{\rightarrow}q$)
and defining
\begin{eqnarray}
N=\frac{1}{2} h + \frac{i \alpha}{4 \ln q}, \quad a=\mu f,\quad 
a^{\dag}=\frac{- i (q - q^{-1})}{\mu (q + q^{-1})} e,
\end{eqnarray}
where $\mu$ is some constant, 
then $N, \ a, \ a^{\dag}$ satisfy the relations (\ref{y}) with $q^2$ in the 
place 
of $q$, and so the $q^2$--boson algebra is isomorphic to a quotient algebra 
of $U_q(sl(2))$.
This realization is similar to that of \cite{yan3}. Finally it should also 
be mentioned that a $q$--deformed differential operator 
algebra was associated with $L$ also possessing a Hopf algebra structure
in \cite{yan2}.

The author of \cite{yan} defines an $R$--matrix of $L$ as the
invertible element that intertwines between the coproduct of (\ref{hop})
and a coproduct $\bar{\Delta}$ that is obtained from that of (\ref{hop})
by changing $q{\rightarrow}q^{-1}$ so that
\begin{eqnarray}
R{\Delta}={\bar{\Delta}}R
\label{yd}
\end{eqnarray}
This definition is claimed to lead to the following $R$--matrix\cite{yan}   
\begin{eqnarray}
R = q^{\left( N - \frac{i \alpha}{\gamma} I \right) \otimes 
\left( N - \frac{i \alpha}{\gamma} I \right) - \frac{1}{2}N{\otimes}N} 
\sum_{k = 0}^{\infty} 
\frac{i^k(1+q^{-1})^{k} q^{-k(k+1)/4}}{\prod_{j = 1}^k \left[ \frac{j}{2} 
\right]_q} (a^{\dag})^k \otimes q^{-kN/2} a^k.
\label{yr}
\end{eqnarray}
Moreover the author states that $R$ satisfies the Yang--Baxter equation
(\ref{dd})(the same $R$ appears also in \cite{yan3}).
Certain comments relative to the above definition of $R$ have to be made.
One would expect that $R$ should be defined by 
relation such as (\ref{da}) as this would justify its nature as an
interwiner between the Hopf algebra structure and the opposite one.
Instead it seems that in \cite{yan} $R$ is not treated as such. 
$\bar{\Delta}$ together with the counit and antipode given in (\ref{hop}) 
does not constitute
a Hopf algebra, not even a coalgebra, and also $\bar{\Delta}$ taken with the
counit and the inverse of S in (\ref{hop}) does not constitute a Hopf 
algebra either. It is the Hopf structure obtained from (\ref{hop}) by setting 
$q{\rightarrow}q^{-1}$ everywhere which is consistent with $\bar{\Delta}$
(obviously this change
leaves (\ref{y}) unaffected). It should be noted that definition
(\ref{yd}) is reminiscent of the case for quantum groups\cite{tol} 
where the definition of the reduced $\bar{R}$--matrix, is given 
by ${\bar{R}}T{\bar{\Delta}}=T{\Delta}{\bar{R}}$
(T being the twist map) and
satisfies relations similar to the Yang--Baxter
equation (the $R$--matrix then can be expressed as a product of $\bar{R}$
with an appropriate $q$--exponentiated function of the Cartan subalgebra 
basis of the Lie algebra and satisfies (\ref{da})). In fact $\bar{R}$ has also
been used by Lusztig\cite{lu}.  
These considerations suggest that (\ref{yd}) alone cannot be used as a 
definition 
either of an $R$--matrix or of a reduced $\bar{R}$--matrix. One can indeed 
verify by direct calculation that $R$ given by (\ref{yr}) {\it does not} 
satisfy any of the relations
 (\ref{yd}), (\ref{dd}), (\ref{da}) or the relation for the reduced 
$\bar{R}$--matrix mentioned above, for example $R{\Delta}(N)=\Delta(N)R$. 
It is the implementation of the 
natural definition of $R$ (\ref{da}), that leads to
the correct $R$--matrix (where $T\Delta$ is compatible with the counit 
and coproduct given in (\ref{hop})) that we shall demonstrate in what
follows using the quantum double construction whose structure
will now be presented. 

Let $A^{\ast}$ denote the dual of $A$ with elements $a^{\ast}$ defined 
by
$(a^{\ast},b) = a^{\ast}(b)$, ${\forall}a{\in}A^{\ast}$ and ${\forall}b 
{\in}A$, where  ( , ) is the natural bilinear form 
$A^{\ast}{\otimes}A{\rightarrow}\IC$ (with $A$ and $A^{\ast}$ regarded
as vector spaces). We assume that
$$A^0 = [ a^{\ast} {\in} A^{\ast} | kera^{\ast} \quad \mbox{contains a 
cofinite two sided
 ideal of A}]$$
is dense in $A^{\ast}$ i.e.
 $(A^0)^\perp = [ a\in A| (b^{\ast} , a) = 0, \quad
{\forall}b^{\ast}{\in} A^0] = (0) $. For $A$ finite dimensional, 
$A^0=A^{\ast}$. Moreover if $A$ is such that the intersection of all
 cofinite two--sided 
ideals is $(0)$ then for every $a{\in} A$ and every $b^{\ast}{\in}
 A^0$ we have
$$a= {\sum}_{s}a_{s}(a_{s}^{\ast}, a)$$
$$b^{\ast} = {\sum}_{s}(b^{\ast}, a_{s})a_{s}^{\ast}$$
where $a_{s}$ and $a_{s}^{\ast}$ are the basis of $A$ and $A^0$ such that
$(a_{s}^{\ast}, a_{t}) = {\delta}_{st}$ and $s = 1, \ 2, ...., \ \dim A$.
Following Sweedler\cite{sweed} we have: 

{\it Theorem}:  $A^0$ becomes a Hopf algebra with multiplication $m^0$, unit
 $u^0$, coproduct $\Delta^0$, antipode $S^0$ and counit $\epsilon^0$ 
defined by:
\begin{eqnarray}
m^0 = {\Delta}^{\ast}|_{A^0{\otimes}A^0}, \quad  
 u^0 = {\epsilon}^{\ast}|_{A^0}\\ 
{\Delta}^0= m^{\ast}|_{A^0}, \quad S^0 = S^{\ast}|_{A^0}, \\  
\epsilon^0(a^{\ast}) = 
(a^{\ast}, 1) \quad {\forall}a^{\ast}{\in}A^0
\end{eqnarray}
where $m^{\ast}$, ${\Delta}^{\ast}$, $ {\epsilon}^{\ast}$, and 
$S^{\ast}$ are the dual maps of $m$, $\Delta$,  $\epsilon$, $S$ 
respectively. The identity element of $A^0$ is given by the counit of $A$.

In what will
follow we shall consider the opposite Hopf algebra structure on $A^0$,
where we keep the same multiplication, unit and counit of $A^0$, but we
use the coproduct $T\Delta^{0}$ and antipode ${(S^0)}^{-1}$
given by
\begin{eqnarray}
(T\Delta^{0}(a^{\ast}), b{\otimes}c ) = (a^{\ast}, c b)\\
({(S^0)}^{-1}(a^{\ast}), b) = (a^{\ast}, S^{-1}(b)), 
\end{eqnarray}
${\forall}a^{\ast}{\in}A^0$ and $b,\  c {\in}A$.
Finally it should be borne in mind that both $A{\otimes}{A^0}$ and 
${A^0}{\otimes}A$ inherit a Hopf algebra structure with respective 
coproducts $\Delta'$ and $\Delta''$ given by
\begin{eqnarray}
\Delta' =
(I{\otimes}\tau{\otimes}I)({\Delta}{\otimes}T\Delta^{0}),\\
\Delta'' =
(I{\otimes}{\tau}^{-1}{\otimes}I)(T\Delta^{0}{\otimes}\Delta),
\end{eqnarray}
where $\tau$ is the twist isomorphism 
$A{\otimes}{A^0}{\rightarrow}{A^0}{\otimes}A$ given by
$\tau(a{\otimes}b^{\ast})= b^{\ast}{\otimes}a$.
 
In the quantum double construction\cite{drin,gr} we first construct the 
vector space
$D(A)$ called the quantum double of $A$ which is the vector space of all 
free products of the form $ab^{\ast}$, 
${\forall}a{\in}A$ and ${\forall}b^{\ast}{\in}A^0$. $D(A)$ is isomorphic 
(as a vector space) with
$A{\otimes}{A^0}$, the isomorphism being given by 
$\psi(a{\otimes}b^{\ast}) = ab^{\ast}$, 
${\forall}a{\in}A$ and ${\forall}b^{\ast}{\in}A^0$. 

$D(A)$ becomes an algebra
by defining all products of the form $b^{\ast}a$, ${\forall}a{\in}A$ 
and ${\forall}b^{\ast}{\in}A^0$, as $b^{\ast}a = \mu(b^{\ast}{\otimes}a)$
where the $\mu$ : ${A^0}{\otimes}A{\rightarrow}D(A)$ given by:
\begin{eqnarray}
A^{0}{\otimes}A\stackrel{(tr{\otimes}I^{{\otimes}2})
((S^0)^{-1}{\otimes}I^{{\otimes}3}){\Delta''}}{\longrightarrow}A^0{\otimes}A
\stackrel{(I^{{\otimes}2}
{\otimes}tr){\Delta''}}{\longrightarrow}A^0{\otimes}A\stackrel{
{\tau}^{-1}}{\longrightarrow}A{\otimes}A^0\stackrel{\psi}
{\rightarrow}D(A)
\end{eqnarray}
where $tr : A^0{\otimes}A {\rightarrow} \IC$ is given by 
$tr(b^{\ast}{\otimes}a) = (b^{\ast}, a)$. Explicitly we have
\begin{eqnarray}
b^{\ast}a= \sum_{(a),(b^{\ast})}({(S^0)}^{-1}((b^{\ast})^{(1)}, a^{(1)})
((b^{\ast})^{(3)}, a^{(3)})a^{(2)}(b^{\ast})^{(2)} 
\end{eqnarray}
Both $A$ and $A^0$ are embedded in 
$D(A)= A{\otimes}{A^0}$ by identifying
$Ia^{\ast}$ and $a{\epsilon}$ with $a^{\ast}$ and $a$
respectively, ${\forall}a^{\ast}{\in}A^{0}$ and ${\forall}a{\in}A$. 

$D(A)$ becomes a quasitriangular Hopf algebra with coproduct ${\Delta}_D$, 
counit ${\epsilon}_D$, and antipode $S_D$ and canonical element $R$ given by
\begin{eqnarray}
{\Delta}_D(ab^{\ast}) =  \Delta(a)(T\Delta^{0})(b^{\ast})\\
{\epsilon}_D(ab^{\ast}) = \epsilon(a){\epsilon}^{0}(b^{\ast})\\
S_D(ab^{\ast}) = {(S^0)}^{-1}(b^{\ast})S(a)\\
R=\sum_{s}a_{s}{\otimes}a^{\ast}_{s} \quad {\in}D(A){\otimes}D(A) 
\label{can}\\
R^{-1} = (S_{D}{\otimes}I)R
\end{eqnarray}
where $[a_{s}]$ and $[a^{\ast}_{t}]$ are bases of $A$ and $A^0$
respectively such that $(a^{\ast}_{s}, a_{t})=\delta_{st}$.

From now on we should always treat the Hopf algebra $L$, as belonging
in the category of quantized universal enveloping algebras. In that
sense we consider $L$ to be spanned by elements of the form
$(N - \frac{i \alpha}{\gamma} I)^{m} . q^{-kN/2}a^k . q^{lN/2}(a^{\dag})^l$ 
($l, m, k = 0, 1, 2, ...$).
We shall denote by $L_{+}$ and $L_{-}$ the Hopf subalgebras
spanned by $(N - \frac{i \alpha}{\gamma} I )^{m} . q^{-kN/2}a^k$ and 
$(N - \frac{i \alpha}{\gamma} I)^{m} . q^{lN/2}(a^{\dag})^l$ respectively.

Following the method just described and similar in spirit with
 \cite{ros} we put $A=L_+$ 
and construct the quantum double $D(L_{+})$ of 
$L_{+}$, to obtained an $R$--matrix for $L$ compatible with
definitions (\ref{da}) and (\ref{dd}). 
 We shall denote by $u$ and $m$ the unit and multiplication 
of $L_+$ while the coproduct $\Delta$, counit $\epsilon$ and antipode $S$ are 
as in (\ref{hop}). Taking $L_+$ to be generated by
$N - (i{\alpha}/{\gamma})I$ and $q^{N/2}a^{\dag}$, we shall denote its
basis by $e_{km} = q^{kN/2} \left( N - \frac{i \alpha}{\gamma} I \right)^m 
(a^{\dag})^k$ with $k,m \geq 0$.
 
As $L_+$ is a coalgebra its dual $L_{+}^{\ast}$ is
necessarily an algebra and via the above theorem, $L_{+}^{0}$ becomes
a Hopf algebra. Let us now construct explicitly $L_{+}^{0}$. As a vector 
space $L_{+}^{0}$ will consist of all linear maps of $L_{+}$ on to
$\IC$.
We take $L_{+}^{0}$ to be generated by the functionals $\nu, \, \beta $ 
on $L_+$ and taking values in $\IC$ defined by 
\begin{eqnarray}
\nu \left( \left( N - \frac{i \alpha}{\gamma} I \right)^n
q^{\frac{lN}{2}} (a^{\dag})^l \right) &=& \frac{\delta_{l0} \delta_{n1} +
i \alpha \delta_{l0} \delta_{n0}}{\gamma},\\
\beta \left( \left( N - \frac{i \alpha}{\gamma} I \right)^n
q^{\frac{lN}{2}} (a^{\dag})^l \right) &=& \frac{e^{-i \alpha/2} \delta_{l1}}
{2^n (1 + q^{-1})}, \nonumber
\label{f}
\end{eqnarray}
and extending by linearity. 

The Hopf algebra structure of $L_{+}^{0}$ can easily be found from
the above theorem. We are interested in the opposite Hopf algebra 
structure of $L_{+}^{0}$
where the multiplication, unit and counit are as in the theorem but
we use the opposite coproduct $T\Delta^{0}$ and antipode $(S^{0})^{-1}$
on $L_{+}^{0}$ given by :
\begin{eqnarray}
T\Delta^{0}(\nu) &=& \nu \otimes 1^{*} + 1^{*} \otimes \nu - 
\frac{i \alpha}{\gamma} 1 \otimes 1^{*},\\
T\Delta^{0}(\beta) &=& \left( \beta \otimes q^{\nu/2} + i q^{-\nu/2} 
\otimes \beta \right) e^{-i \alpha / 2},\\
(S^{0})^{-1}({\nu}) &=& -{\nu} + (2i{\alpha}/{\gamma}) 1^{*},\\
(S^{0})^{-1}({\beta}) &=& - \frac{1}{q^{1/2}}{\beta}
\end{eqnarray}
where $1^{*}$ is the identity on $L_{+}^{0}$
(i.e. $u^{0}(1)=1^{*}$). Moreover as $L_{+}^{0}$ inherits a Lie algebra 
structure, with the non--zero Lie bracket given by
\begin{eqnarray} 
[\nu,\beta] &=& -\beta,
\end{eqnarray}
as can be seen using (\ref{f}).
It is convenient to define a basis of $L_{+}^{0}$ given by
\begin{eqnarray}
e_{km}^{\ast}&=& 
\left( \nu - \frac{i \alpha}{\gamma} 1^{*} \right)^{m} q^{-k\nu/2}\beta^k,
\quad \mbox{so that}\\ 
e_{km}^{\ast}(e_{ln}) &=& \left( q^{-k\nu/2} \left( \nu - \frac{i \alpha}{\gamma} 1^{*} \right)^m \beta^k 
\right) \left( q^{lN/2} \left( N - \frac{i \alpha}{\gamma} I \right)^n 
(a^{\dag})^{l} \right) \nonumber\\
&=& \delta_{kl} \delta_{mn} \frac{n! (-i)^k q^{k(k+1)/4}}{\gamma^n} . 
\prod_{j = 1}^k \left[ \frac{j}{2} \right]_q.
\label{n}
\end{eqnarray}
Observe that the map $\nu{\rightarrow}N$ and 
$\beta{\rightarrow}a$ defines an isomorphism $L_{+}^{0} \cong L_{-}$.

Considering now the quantum double $D(L_{+})=L_{+}{\otimes}L_{+}^{0}$.
It is the vector space spanned by all free products $ab^{\ast}$ which
becomes an algebra by defining all products of the form $b^{\ast}a$ as
was indicated above. Moreover as was stated it is a quasitriangular
Hopf algebra and thus by appropriately normalising the elements 
$e_{km}$ and $e_{lm}^{\ast}$ using (\ref{n}), the canonical element $R$ given in (\ref{can})
is realised as
\begin{eqnarray}
R=q^{\left( N - \frac{i \alpha}{\gamma} I \right) \otimes 
\left( \nu - \frac{i \alpha}{\gamma} 1^{*} \right)} \sum_{k = 0}^{\infty} 
\frac{i^k q^{-k(k+1)/4}}{\prod_{j = 1}^k \left[ \frac{j}{2} \right]_q} 
q^{kN/2} (a^{\dag})^k \otimes q^{-k\nu/2} \beta^k \in D(L_{+}).
\end{eqnarray}
The relation of our $q$--boson algebra $L$ with $D(L_{+})$ can now be
obtained. Observe firstly, that we have the
following quantum double intertwining relations between $L_{+}$ and
$L_{+}^{0}$: 
\begin{eqnarray}
{[}N,{\nu}{]} &=& 0,\\
{[}N,{\beta}{]} &=& - {\beta},\\
{[}{\nu},a^{\dag}{]} &=& a^{\dag},\\
{[}{\beta},a^{\dag}{]} &=& [ \frac{{\nu} + N + 1}{2}]_q - [\frac{{\nu} + N - 1}{2} ]_q.  
\end{eqnarray}
Secondly note that the element $\nu - N$ is central in $D(L_{+})$ and 
generates a two--sided Hopf ideal, call it $M$.
The quotient Hopf algebra $D(L_{+})/M$, in which $\nu = N$
can be identified with the $q$--boson algebra $L$ by identifying $\beta$
with $a$ and $1^{*}$ with $I$. The $R$--matrix is now given by
\begin{eqnarray}
R = q^{\left( N - \frac{i \alpha}{\gamma} I \right) \otimes 
\left( N - \frac{i \alpha}{\gamma} I \right)} \sum_{k = 0}^{\infty} 
\frac{i^k q^{-k(k+1)/4}}{\prod_{j = 1}^k \left[ \frac{j}{2} \right]_q} 
q^{kN/2} (a^{\dag})^k \otimes q^{-kN/2} a^k.
\label{rr}
\end{eqnarray}
and satisfies the Yang--Baxter equation. 

The differences between (\ref{yr}) and (\ref{rr}) can now be read off
and the cause of the inconsistency of (\ref{yr}) is obvious. 
The representation theory of (\ref{y}) in the spirit of \cite{sg},
the relation
of (\ref{rr}) with the $R$--matrix of $U_{q}(sl(2))$, its connection 
with representations of the braid group and possible relation with
link invariances are under investigation. Finally we conclude with the 
observation that besides the Hopf 
algebra structure given by (\ref{hop}) a more general one exists 
given by 
\begin{eqnarray}
\Delta(N) = N \otimes I + I \otimes N + \beta I \otimes I, \nonumber\\
\Delta(a) = \left( a \otimes q^{m N} \pm (-1)^K i q^{(m \pm 1) N} \otimes a 
\right) e^{i \pi (2K+1) m / 2},\nonumber\\
\Delta(a^{\dag}) = \left( a^{\dag} \otimes q^{- (m \pm 1) N} \pm (-1)^K i 
q^{- m N} \otimes a^{\dag} \right) e^{i \pi (2K+1) (m \pm 1) / 2},
\nonumber\\
\varepsilon(N) = -{\beta}, \quad \varepsilon(a) 
= \varepsilon(a^{\dag}) = 0,\quad \varepsilon(I)=1, \nonumber\\
S(N) = - N -2{\beta}I, \nonumber\\
S(a) = {\pm} i(-1)^{K}q^{-(m{\pm}1)N} aq^{mS(N)},\nonumber\\
S(a^{\dag}) = {\pm} i(-1)^{K}q^{mN} a^{\dag}q^{-(m{\pm}1)S(N)}. 
\label{gen1}
\end{eqnarray}
where $\beta = \frac{i \pi (2K+1)}{2 \gamma}$, $\gamma=\ln q$, $m$ is an
integer or half--integer, and $K$ is an integer. The corresponding 
$R$--matrix is now given by 
\begin{eqnarray}
R = q^{\mp \left( N + \beta I \right) \otimes \left( N + \beta I 
\right)} \sum_{k = 0}^{\infty} 
\frac{(\pm i)^k (-1)^{Kk}}{\prod_{j = 1}^k \left[ \frac{j}{2} \right]_q} 
q^{-m k^2 \mp \frac{1}{4} k(k-1)} q^{mkN} (a^{\dag})^k \otimes q^{-mkN} a^k.
\label{gen2}
\end{eqnarray}
The choices $m=1/2$, $K=-2{\kappa}-1$ $(\kappa{\in}\IZ)$ and the lower
sign in (\ref{gen1}) and (\ref{gen2}), lead directly  to (\ref{hop}) and 
(\ref{rr}) respectively.

The authors would like to thank P. D. Jarvis, A. J. Bracken 
for useful comments and support during the completion of this
letter. One of us (I.T) would also like to thank R. Zhang and A. Ram 
for fruitful discussions during the ``Conference on Lie Theory", 
27 Nov- 1 Dec 1995, Institute for Theoretical Physics, Adelaide, 
where the content of this letter was reported.


\begin{thebibliography}{99}
\bibitem{drin}Drinfeld V G 1986 `Quantum Groups' in Proc. ICM Berkeley,
{\bf1} 798
\bibitem{j}Jimbo M. 1985 Lett. Math. Phys. {\bf10} 63
1986 Lett. Math. Phys. {\bf11} 247; 1987 Commun. Math. Phys. {\bf102} 537
\bibitem{q}`Quantum Groups' 1990 Proceedings of the Argonne     Workshop
World Scientific 1990 ed. Curtright T, Fairlie D and Zachos C.
\bibitem{kr}Kulish P P and Reshetikhin N Yu 1989 Lett. Math. Phys. {\bf18} 143
\bibitem{tol}Khoroshkin S M and Tolstoy V N 1991 Comm. Math. Phys.
{\bf141} 599
\bibitem{gr1}Bracken A J, Gould M D and Zhang R B 1990 Mod. Phys. Lett. 
{\bf A5} 831; Bracken A J, Gould M D and Tsohantjis I 1993 J.
Math. Phys. {\bf34} 1654
\bibitem{sweed}Sweedler M E 1969 Hopf Algebras (Benjamin, NY)
\bibitem{kur}Kuryshkin W 1980 Ann. Found. L de Broglie {\bf5} 111
\bibitem{j1}Jannussis A et al 1982 Hadronic J. {\bf5} 1923
Jannussis A, Brodimas G, Sourlas D and Zisis, V 1981 Lett. Nuovo Cimento
{\bf30} 123;Jannussis A 1991 Hadronic J. {\bf14} 257
\bibitem{mac}Macfarlane A J 1989 J. Phys. A: Math. and Gen. {\bf22} 4581 
\bibitem{bid}Biedenharn L C 1989 J. Phys. A: Math. Gen. {\bf22} L873
\bibitem{bid2}Sun C-P and Fu H-C 1989 J. Phys. A: Math. Gen. {\bf22} L983;
Chaichian M, Kulish P P and Lukierski J 1990 Phys. Lett. {\bf237B} 401
\bibitem{ch}Chaichian M and Kulish P P 1990 Phys. Lett. {\bf234B} 72
\bibitem{che}Chaichian M and Ellinas D 1990 J. Phys. A: Math. and Gen.
{\bf23} L291
\bibitem{dav}Bracken A J, McAnally D S, Zhang R B and Gould M D 1991 
J. Phys. A: Math. Gen. {\bf24} 1379
\bibitem{yan2}Hong Yan 1991 J. Phys. A:Math. Gen.{\bf24} L409 
\bibitem{yan}Hong Yan 1990 J. Phys. A:Math. Gen. {\bf23} L1155
\bibitem{oh}Oh C H and Singh K 1994 J. Phys. A: Math. Gen. {\bf27} 5907
\bibitem{ros}Rosso M 1989 Commun. Math. Phys. {\bf124} 307
\bibitem{gr}Gould M D, Zhang R B and Bracken A J 1993 Bull. Austral.
Math. Soc {\bf47} 353; Gould M D 1993 Bull. Austral. Math. Soc {\bf48} 275;
Tsohantjis I and Gould M D 1994 Bull. Austral. Math. Soc {\bf49} 177
\bibitem{ng}Ng Y T 1990 J. Phys. A: Math. and Gen.{\bf23} 1023
\bibitem{b}Brodimas G `Lie Admissible $Q$-Algebras and Quantum Groups' 
Phd Thesis 1991 University of Patras, Patras, Greece
\bibitem{da}Daskaloyannis C 1991 J. Phys. A: Math. Gen. {\bf24} L789
\bibitem{jo}J{\o}rgensen P E T, Schmitt L M and Werner R F 1994 Pacific J. 
Math. {\bf 165} 131; 1995 J. Fun.Anal. {\bf134} 33
\bibitem{cgp}Chaichian M, Grosse H and Pre\v{s}najder P 1994 
J. Phys. A: Math. Gen. {\bf27} 2045; Chaichian M, Felipe Gonzalez R and 
Pre\v{s}najder P 1995 J. Phys. A: Math. Gen. {\bf28} 2247
\bibitem{yan3}Hong Yan 1991 Phys. Lett {\bf262B} 459
\bibitem{kd}Kulish P P and Damaskinsky E V 1990 J. Phys. A: Math. Gen. 
{\bf23} L415
\bibitem{qexp}McAnally D S 1995 J. Math. Phys. {\bf 36} 546
\bibitem{qgam}McAnally D S 1995 J. Math. Phys. {\bf 36} 574
\bibitem{lu}Lusztig G 1993 ``Introduction to Quantum Groups" 
(Boston, Birkhauser)
\bibitem{sg}Sun Chang-Pu and Ge Mo-Lin (1991) J. Math. Phys. {\bf32} 595




\end{thebibliography}
\end{document}